\documentclass[11pt,draft,amsmath,amscd,amsbsy,amssymb]{amsart}
 \relax

\chardef\bslash=`\\ 

\makeatletter \def\verbatim{\interlinepenalty\@M \@verbatim
\leftskip\@totalleftmargin\advance\leftskip2pc
\frenchspacing\@vobeyspaces \@xverbatim} \makeatother \hfuzz1pc

\makeatletter \def\dgt@k{\dg@DX=-3 \dg@DY=2 \dg@SIZE=3}
\makeatother

\makeatletter \def\dgt@kk{\dg@DX=3 \dg@DY=-1 \dg@SIZE=3}

\theoremstyle{plain} \newtheorem{thm}{Theorem}[section]
\newtheorem{cor}[thm]{Corollary}
\newtheorem{lemma}[thm]{Lemma}
\newtheorem{prop}[thm]{Proposition}

\theoremstyle{definition} \newtheorem{rem}[thm]{Remark}
\newtheorem{defin}[thm]{Definition} \newtheorem{que}[thm]{Question}

 \newtheorem{exam}[thm]{Example}

\newcommand{\R}{\mathbb R_{\mathrm{max}}}

\newcommand{\comp}{\operatorname{\mathbf{Comp}}}

\newcommand{\diam}{\operatorname{\mathrm{diam}}}

\newcommand{\pr}{\mathrm{pr}}
\newcommand{\id}{\mathrm{id}}

\usepackage[all]{xy}

\begin{document}

\title[Idempotent probability measures]
{Idempotent probability measures, I}
\author{M. Zarichnyi}
\address{Department of Mechanics and Mathematics, Lviv National University,  Universtytetska Str. 1,
79000 Lviv, Ukraine} \email{mzar@litech.lviv.ua}

\thanks{}
\subjclass[2000]{Primary 54C65, 52A30; Secondary 28A33}



\begin{abstract} The set of all idempotent
probability measures (Maslov measures) on a compact Hausdorff
space endowed with the weak* topology determines is functorial on
the category $\comp$ of compact Hausdorff spaces. We prove that
the obtained functor is normal in the sense of E. Shchepin. Also,
this functor is the functorial part of a monad on $\comp$. We
prove that the idempotent probability measure monad contains the
hyperspace monad as its submonad. A counterpart of the notion of
Milyutin map is defined for the idempotent probability measures.
Using the fact of existence of Milyutin maps we prove that the
functor of idempotent probability measures preserves the class of
open surjective maps. Unlikely to the case of probability
measures, the correspondence assigning to every pair of idempotent
probability measures on the factors the set of measures on the
product with these marginals, is not open.
\end{abstract}

\maketitle
\section{Introduction}
In this paper we  investigate the functor of idempotent
probability measures in the category of compact Hausdorff spaces.
Hopefully, this functor will play in the idempotent analysis a
role similar to that of the probability measure functor in the
classical functional analysis.

The notion of idempotent (Maslov) measure finds important
applications in different part of mathematics, mathematical
physics and economics (see the survey article \cite{L} and the
bibliography therein). In particular, these measures arise in
dynamical optimization \cite{be}. An analogy between Maslov
integration and optimization is indicated in \cite{ad}. It is
noted in \cite{au} that ``the use of Maslov measures for
encapsulating some aspects of uncertainty can be as relevant to
most economic problem than the use of classical probability
theory''.

According to an idempotent correspondence principle \cite{L},
``there exists a heuristic correspondence between important,
interesting, and useful constructions and results of the
traditional mathematics over fields and analogous constructions
and results over idempotent semirings and semifields (i.e.,
semirings and semifields with idempotent addition)''.

E. Shchepin \cite{S} introduced the class of normal functors
acting in the category $\comp$ of compact Hausdorff spaces and
continuous maps and showed that the probability measure functor
$P$ is normal (see \cite{F} for a detailed proof of this fact).
The aim of this paper is to establish the property of normality
for the functor of idempotent probability measures. Moreover, we
show that the latter functor determines a monad (see the
definition below) in the category $\comp$. It is also proved that
the functor of idempotent probability measures contains the
hyperspace functor as a subfunctor. Moreover, the hyperspace monad
can be embedded as a submonad in the monad of idempotent
probability measures. The proofs of these results are based on
embeddings of the spaces of idempotent probability measures in the
spaces of order-preserving functionals defined by Radul \cite{R}.
The monad structure for the functor of idempotent probability
measures is tightly connected with the max-plus convex compact
subsets in euclidean spaces as these sets endowed with the
idempotent barycenter maps turn out to be algebras for the monad
of idempotent probability measures.

A Milyutin map is a map of topological spaces that admits an
averaging operator. Equivalently, the Milyutin maps are precisely
those admitting the probability-measure-valued selection. We also
define a natural counterpart for the notion of Milyutin map for
the spaces of idempotent probability measures. Using this notion
we prove that the functor of idempotent probability measures is
open, i.e. it preserves the class of open surjective maps. Note
that the openness of the functor of probability measures is
established by Ditor and Eifler \cite{DE}.

Eifler \cite{E}  proved that the correspondence that assigns to
every pair of probability measures on factors the set of
probability measures with given marginals, is open (see also
\cite{Z}). We prove that this is not the case for the  idempotent
probability measures.

\section{Preliminaries}\label{s:1}

Let $X$ be a compact Hausdorff space.  By $C(X)$ we denote the
Banach space of continuous functions on $X$ endowed with the
$\sup$-norm. For any $c\in\mathbb R$ we denote by $c_X$ the
constant function on $X$ taking the value $c$. By $w(X)$ we denote
the weight of a topological space $X$.

In the sequel, by functor we mean a covariant functor. The
probability measure functor acting in the category $\comp$ is
denoted by $P$. See, e.g. \cite{F} for the properties of the
functor $P$.

By $\exp$ we denote the hyperspace functor acting in the category
$\comp$. Given a compact Hausdorff space $X$, the space $\exp X$
is defined as the set of all nonempty closed subsets in $X$
endowed with the Vietoris topology. A base of this topology is
formed by the sets of the form $$\langle
U_1,\dots,U_n\rangle=\{A\in \exp X\mid A\subset\cup_{i=1}^nU_i,\
A\cap U_i\neq\emptyset,\ i=1,\dots,n\}.$$ If $f\colon X\to Y$ is a
continuous map, then $\exp f\colon \exp X\to\exp Y$ is defined by
$\exp f(A)=f(A)$, $A\in\exp X$.

Let $\R=\mathbb R\cup\{-\infty\}$ endowed with the metric
$\varrho$ defined by $\varrho(x,y)=|e^x-e^y|$. Let also
$\R^n=(\R)^n$.

A functor in the category $\comp$ is called {\em open} if it
preserves the class of open surjective maps. Recall that a map of
topological spaces is called {\em open} if the image of every open
set is open. For a surjective map $f\colon X\to Y$ of compact
Hausdorff spaces the openness of $f$ is equivalent to the
continuity of the map $y\mapsto f^{-1}(y)\colon Y\to \exp X$. If
moreover $X$ and $Y$ are metrizable, then $f$ is open if and only
if, for any sequence $(y_i)_{i=1}^\infty$ converging to $y\in Y$
and every $x\in X$ such that $f(x)=y$, there exists a sequence
$(x_i)_{i=1}^\infty$ converging to $x$ and such that $f(x_i)=y_i$,
for every $i$.

Following the style of idempotent mathematics (see, e.g.,
\cite{L,LMS,MS}) we denote by $\odot\colon \mathbb R\times C(X)\to
C(X)$ the map acting by $(\lambda,\varphi)\mapsto
\lambda_X+\varphi$, and by $\oplus\colon C(X)\times C(X)\to C(X)$
the map acting by $(\varphi,\psi)\mapsto \max\{\varphi,\psi\}$.

For each $c\in\mathbb R$ by $c_X$ we denote the constant function
from $C(X)$ defined by the formula $c_X(x)=c$ for each $x\in X$.

\begin{defin}\label{d:1} A functional $\mu\colon C(X)\to\mathbb R$ is called an {\em
idempotent probability measure} (a {\em Maslov measure}) if
\begin{enumerate}
\item $\mu(c_X)=c$;
\item $\mu(c\odot\varphi)=c\odot\mu(\varphi)$;
\item $\mu(\varphi\oplus\psi)=\mu(\varphi)\oplus\mu(\psi),$
\end{enumerate}
for every $\varphi,\psi\in C(X)$.
\end{defin}

The number $\mu(\varphi)$ is the {\em Maslov integral} of
$\varphi\in C(X)$ with respect to $\mu$.

Let $I(X)$ denote the set of all idempotent probability measures
on $X$. We endow $I(X)$ with the weak* topology. A base of this
topology is formed by the sets
$$\langle\mu; \varphi_1,\dots,\varphi_n; \varepsilon\rangle=\{\nu\in I(X)\mid
|\mu(\varphi_i)-\nu(\varphi_i)|<\varepsilon,\ i=1,\dots,n\},$$
where $\mu\in I(X)$, $\varphi_i\in C(X)$, $i=1,\dots,n$,  and
$\varepsilon>0$.

The following is an example of an idempotent probability measure.
Let $x_1,\dots,x_n\in X$ and $\lambda_1,\dots,\lambda_n\in\R$ be
numbers such that $\max\{\lambda_1,\dots,\lambda_n\}=0$. Define
$\mu\colon C(X)\to\mathbb R$ as follows:
$\mu(\varphi)=\max\{\varphi(x_i)+\lambda_i\mid i=1,\dots,n\}$.  As
usual, for every $x\in X$, we denote by $\delta_x$ (or
$\delta(x)$) the functional on $C(X)$ defined as follows:
$\delta_x(\varphi)=\varphi(x)$, $\varphi\in C(X)$ (the Dirac
probability measure concentrated at $x$). Then one can write
$\mu=\oplus_{i=1}^n\lambda_i\odot\delta_{x_i}$.

In order to establish properties of the set of idempotent
probability measures, we use those of the order-preserving
functionals. In \cite{R}, T. Radul  considered the set of
order-preserving functionals  on  compact Hausdorff spaces.

A functional (which is not supposed a priori to be either linear
or continuous) $\nu\colon C(X)\to\mathbb R$ is called
\begin{enumerate}
\item
{\it weakly additive\/} if for each $c\in\mathbb R$ and
$\varphi\in C(X)$ we have $\nu(\varphi+c_X)=\nu(\varphi)+c$;
\item {\it
order-preserving\/} if for each $\varphi, \psi\in C(X)$ with
$\varphi\le\psi$ we have $\nu(\varphi)\le\nu(\psi)$;
\item  {\it normed} if $\nu(1_X)=1$.

\end{enumerate}

The space of real numbers $\mathbb R$ is endowed with the~standard
metric. The following fact is established in \cite{R}.
\begin{lemma}
Each order-preserving weakly additive functional is a
non-expanding map.
\end{lemma}

For a compact Hausdorff space $X$, we denote by $O(X)$ the set of
all order-preserving weakly additive normed functionals in $C(X)$.
It is easy to see that for each $\nu\in O(X)$ and $c\in\mathbb R$
we have $\nu(c_X)=c$. Therefore, $I(X)\subset O(X)$, for every
compact Hausdorff space $X$.

\begin{prop} The set $I(X)$ is closed in $O(X)$.
\end{prop}
\begin{proof} Suppose that $\mu\in O(X)\setminus I(X)$. Then there
exist $\varphi,\psi\in C(X)$ such that
$a=\mu(\varphi\oplus\psi)>\mu(\varphi)\oplus\mu(\psi)=b$. Then
$\mu\in \left\langle\mu;\varphi\oplus\psi,\varphi,\psi;
\frac{a-b}{2}\right\rangle\subset O(X)\setminus I(X)$. We see that
the complement of $I(X)$ is an open set in $O(X)$.

\end{proof}

Since $O(X)$ is known to be compact Hausdorff, we conclude that so
is the space $I(X)$.

Given a map $f\colon X\to Y$ of compact Hausdorff spaces, the map
$O(f)\colon O(X)\to O(Y)$ by the formula
$O(f)(\mu)(\varphi)=\mu(\varphi f)$, for every $\varphi\in C(Y)$.

\begin{prop} Let $f\colon X\to Y$ be a continuous map of compact
Hausdorff spaces. Then $O(f)(I(X))\subset I(Y)$.
\end{prop}
\begin{proof} Let $\mu\in I(X)$ and
$\varphi,\psi\in C(Y)$. Clearly, $O(f)(\mu)(c_X)=c$,  $c\in\mathbb
R$, and $O(f)(\mu)(\lambda\oplus \varphi)=\lambda\oplus
O(f)(\mu)(\varphi)$, $\lambda\in \mathbb R$.  We have also
\begin{align*} O(f)(\mu)(\varphi\oplus\psi)=&\mu((\varphi\oplus\psi)f)=\mu((f\varphi\oplus
f\psi))\\= &\mu(f\varphi)\oplus
\mu(f\psi)=O(f)(\mu)(\varphi)\oplus O(f)(\mu)(\psi).\end{align*}
Thus, $O(f)(\mu)\in I(Y)$.

\end{proof}

We denote by $I(f)\colon I(X)\to I(Y)$ the restriction map $O(f)|
I(X)\colon I(X)\to I(Y)$. Note that, if
$\mu=\oplus_{i=1}^n\lambda_i\odot\delta_{x_i}\in I(X)$, then
$I(f)(\mu)=\oplus_{i=1}^n\lambda_i\odot\delta_{f(x_i)}\in I(Y)$.

It is evident that the construction $I$ determines a covariant
functor in the category $\comp$.

\begin{prop}\label{p:emb} The functor $I$ preserves the class of embeddings.
\end{prop}
\begin{proof} This directly follows from the fact that the functor
$O$ preserves the embeddings \cite{R}.
\end{proof}

As usual, given a closed subset $A$ of a compact Hausdorff space
$X$, we identify $I(A)$ with the subspace $I(i)(I(A))$ of $I(X)$,
where by $i\colon A\to X$ we denote the inclusion map.

The proof of the following statement involves  the max-plus
version of the Hahn-Banach theorem \cite{LMS}. For the sake of
completeness, we provide an alternative proof of its special
version.

We say that a subset $L$ of $C(X)$ is a {\em max-plus} linear
subspace of $C(X)$ if
\begin{enumerate}
\item $c_X\in L$ for every $c\in\mathbb R$;
\item $\lambda\odot\varphi\in L$, for every $\lambda\in\mathbb R$
and $\varphi\in L$;
\item $\varphi\oplus\psi\in L$, for every $\varphi,\psi\in L$.
\end{enumerate}

\begin{lemma}\label{l:hb0} Let $L$ be a max-plus linear subspace of $C(X)$. Let $\mu\colon L\to\mathbb R$ be a functional that
satisfies the conditions of Definition \ref{d:1} (with $C(X)$
replaced by $L$). For any $\varphi_0\in C(X)\setminus L$, there
exists an extension of $\mu$ onto the minimal linear max-plus
subspace $L'$ containing $L\cup\{\varphi_0\}$ that satisfies the
conditions of Definition \ref{d:1}.
\end{lemma}
\begin{proof} For any $\varphi\in L'$, define $\mu(\varphi)=\inf\{\mu(\psi)\mid\psi\in L,\
\varphi\le\psi\}$. It is clear that $\mu$ is well-defined.

Given $\varphi\in L'$ and $\lambda\in\mathbb R$, we see that
\begin{align*}\mu(\lambda\odot\varphi)=&\inf\{\mu(\psi)\mid\psi\in L,\
\lambda\odot\varphi\le\psi\}\\=&\inf\{\mu(\lambda\odot\psi')\mid\psi'\in
L,\ \lambda\odot\varphi\le\lambda\odot\psi'\}\\=&
\lambda\odot\inf\{\mu(\psi')\mid\psi'\in L,\ \varphi\le\psi'\}\\=&
\lambda\odot\mu(\varphi).\end{align*}

Note also that $\mu(\varphi_1)\le\mu(\varphi_2)$, whenever
$\varphi_1\le\varphi_2$, $\varphi_1,\varphi_2\in L'$.

Now, given $\varphi_1,\varphi_2\in L'$, we see that
\begin{align*}\mu(\varphi_1)\oplus\mu(\varphi_2)=&\inf\{\mu(\psi)\mid\psi\in L,\
\varphi_1\le\psi\}\oplus\inf\{\mu(\psi')\mid\psi'\in L,\
\varphi_2\le\psi'\}\\ =&
\inf\{\mu(\psi)\oplus\mu(\psi')\mid\psi,\psi'\in L,\
\varphi_1\le\psi,\ \varphi_2\le\psi'\}\\ \ge&
\inf\{\mu(\psi\oplus\psi')\mid\psi,\psi'\in L,\
\varphi_1\oplus\varphi_2\le\psi\oplus\psi'\}\\=&
\mu(\varphi_1\oplus\varphi_2).\end{align*} On the other hand,
since $\varphi_i\le \varphi_1\oplus\varphi_2$, $i=1,2$, from the
above remark it follows that $\mu(\varphi_i)\le
\mu(\varphi_1\oplus\varphi_2)$, $i=1,2$, and therefore
$\mu(\varphi_1)\oplus\mu(\varphi_2)=\mu(\varphi_1\oplus\varphi_2)$.
\end{proof}
\begin{lemma}\label{l:hb} Let $L$ be a max-plus linear subspace of $C(X)$.
Let $\mu\colon L\to\mathbb R$ be a functional that satisfies the
conditions of Definition \ref{d:1} (with $C(X)$ replaced by $L$).
Then there exists $\nu\in I(X)$ such that $\nu|L=\mu$.
\end{lemma}
\begin{proof}
We  apply Lemma \ref{l:hb0} and the Kuratowski-Zorn lemma to
obtain a maximal extension $\mu'$ of $\mu$ onto a max-plus linear
subspace $L'$ of $C(X)$. If $L'\neq C(X)$ and $\varphi_0\in
C(X)\setminus L'$, then we can extend $\mu'$ onto the minimal
max-plus linear subspace of $C(X)$ containing
$L'\cup\{\varphi_0\}$, which contradicts the maximality.
\end{proof}

\begin{prop}\label{p:onto}  The functor $I$ preserves the class of the onto maps.
\end{prop}
\begin{proof} Let $f\colon X\to Y$ be an onto map of compact
Hausdorff spaces. Let $C(f)=\{\varphi f\mid \varphi\in C(Y)\}$.
Given $\mu\in I(X)$, we consider the map $\mu'\colon C(f)\to
\mathbb R$ defined by the formula $\mu'(\varphi f)=\mu(\varphi)$.
Obviously, $\mu$ satisfies the conditions from Definition
\ref{d:1} (with $C(X)$ replaced by $C(f)$). By Lemma \ref{l:hb},
there exists an extension $\nu$ of $\mu'$ onto $C(X)$ satisfying
the conditions from Definition \ref{d:1}. We then have
$\mu=I(f)(\nu)$.
\end{proof}

\begin{prop}\label{p:int} The functor $I$ preserves the intersections in the
sense that $I(A\cap B)=I(A)\cap I(B)$, for any closed subsets
$A,B$ of a compact Hausdorff space $X$.
\end{prop}
\begin{proof}  Since the functor $O$
preserves intersections \cite{R}, we see that
\begin{align*}I(A\cap B)=& O(A\cap B)\cap I(A\cap B)=(O(A)\cap
O(B))\cap I(A\cap B)\\= &(O(A)\cap I(A\cap B))\cap (O(B)\cap
I(A\cap B))\\=&I(A)\cap I(B).\end{align*}
\end{proof}

If a functor $F$ in $\comp$ preserves the class of embeddings and
intersections, one can define the notion of support for it.
Namely, given $a\in F(X)$, we call the set
$$\mathrm{supp}(a)=\cap\{Y\mid Y\text{ is a closed subset of }X\text{ and }a\in F(Y) \}$$
the {\em support} of $a$. Note that if
$\mu=\oplus_{i=1}^n\lambda_i\odot\delta_{x_i}$, then
$\mathrm{supp}(\mu)=\{x_i\mid\lambda_i>-\infty\}$.

Given a closed subset $A$ of $X$ and $\mu\in I(X)$, we have:
$\mathrm{supp}(\mu)\subset A$ if and only if, for any two
functions $\varphi,\psi\in C(X)$ with $\varphi|A=\psi|A$, we have
$\mu(\varphi)=\mu(\psi)$.

Let $\mathcal S=\{X_\alpha,p_{\alpha\beta};\mathcal A\}$ be an
inverse system over a directed set $\mathcal A$. For any
$\alpha\in\mathcal A$, let $p_\alpha\colon X=\varprojlim\mathcal
S\to X_\alpha$ denote the limit projection. By $I(\mathcal S)$ we
denote the inverse system
$\{I(X_\alpha),I(p_{\alpha\beta});\mathcal A\}$.

The following property is a counterpart of the Kolmogorov theorem
for probability measures.

\begin{prop}\label{p:conti} The map $h=(I(p_\alpha))_{\alpha\in\mathcal A}\colon
I(X)\to\varprojlim I(\mathcal S)$ is a homeomorphism.
\end{prop}

\begin{proof} It easily follows from the results of Radul \cite{R}
that the map $h$ is an embedding. We are going to show that $h$ is
an onto map. Let $(\mu_\alpha)_{\alpha\in\mathcal A}\in
\varprojlim I(\mathcal S)$. By \cite[Proposition 4]{R}, there
exists $\mu\in O(X)$ such that $O(p_\alpha(\mu))=\mu_\alpha$, for
any $\alpha\in\mathcal A$.  Let $C'=\{\varphi p_\alpha\mid
\varphi\in C(X_\alpha),\ \alpha\in\mathcal A\}$. Given
$\varphi,\psi\in C'$, one can write $\varphi=\varphi' p_\alpha$,
$\psi=\psi' p_\alpha$, for some $\alpha\in\mathcal A$, whence
\begin{align*} \mu(\varphi\oplus\psi)= &\mu((\varphi' p_\alpha)\oplus (\psi'
p_\alpha))=
O(p_\alpha)(\mu)(\varphi'\oplus\psi')\\=&\mu_\alpha(\varphi'\oplus\psi')=
\mu_\alpha(\varphi')\oplus\mu_\alpha(\psi')\\ =&\mu(\varphi'
p_\alpha)\oplus \mu(\psi'
p_\alpha)=\mu(\varphi)\oplus\mu(\psi).\end{align*}

Since, by the Stone-Weierstrass theorem, the set $C'$ is dense in
$C(X)$ and the operation $\oplus$ is continuous, we conclude that
$\mu(\varphi\oplus\psi)=\mu(\varphi)\oplus\mu(\psi)$ for all
$\varphi,\psi\in C(X)$. One can similarly prove that
$\mu(\lambda\odot\varphi)=\lambda\odot\mu(\varphi)$, for all
$\varphi\in C(X)$ and $\lambda\in \mathbb R$.

Thus, $\mu\in I(X)$ is as required.
\end{proof}

E. Shchepin \cite{S} calls the just established property of the
functor $I$ the {\em continuity} of $I$.

Suppose that a functor $F$ in the category $\comp$ preserves the
class of embeddings. We say that $F$ {\em preserves preimages} if,
for every map $f\colon X\to Y$ and every closed subset $B$ of $Y$,
we have $F(f^{-1}(B)=(F(f))^{-1}(F(B))$.

\begin{prop}\label{p:prei} The functor $I$ preserves preimages.
\end{prop}

\begin{proof} Assume the contrary and let $f\colon X\to Y$ be a morphism in $\comp$, $\mu\in I(X)$,
$B$ be a closed subset in $Y$ such that $I(f)(\mu)\in I(B)$ while
$\mu\notin I(f^{-1}(B))$. There exist $\varphi,\psi\in C(X)$ such
that $\varphi|f^{-1}(B)=\psi|f^{-1}(B)$ and
$\mu(\varphi)\neq\mu(\psi)$. Let $c=|\mu(\varphi)-\mu(\psi)|$.

There exists a neighborhood $U$ of $B$ in $Y$ such that
$\|\varphi|f^{-1}(\bar U)-\psi|f^{-1}(\bar U)\|<(c/3)$. There
exist functions $\varphi_1,\psi_1\in C(X)$ satisfying the
properties:
\begin{enumerate}
\item $\varphi_1|f^{-1}(\bar U)=\varphi|f^{-1}(\bar U)$, $\psi_1|f^{-1}(\bar U)=\psi|f^{-1}(\bar U)$;
\item $\varphi_1\le\varphi$, $\psi_1\le\psi$;
\item $\|\varphi_1-\psi_1\|<(c/3)$.

 One can easily demonstrate that there exist functions
$\varphi_2,\psi_2\in C(X)$ satisfying the properties:

\item $\varphi_2|(X\setminus f^{-1}(\bar U))=\varphi|(X\setminus f^{-1}(\bar U))$,
$\psi_1|(X\setminus f^{-1}(\bar U))=\psi|(X\setminus f^{-1}(\bar
U))$;
\item $\varphi_2\le\varphi$, $\psi_2\le\psi$;
\item $\varphi_2|f^{-1}(B)=\psi_2|f^{-1}(B)=\min\{\inf\varphi,\inf\psi\}-1$.

Since
$$\mu(\varphi)=\mu(\varphi_1\oplus\varphi_2)=\mu(\varphi_1)\oplus\mu(\varphi_2),\
\mu(\psi)=\mu(\psi_1\oplus\psi_2)=\mu(\psi_1)\oplus\mu(\psi_2),$$
from property (3) and from the choice of $\varphi,\psi$ it follows
that $\mu(\varphi_2)\neq \mu(\psi_2)$.

 There exist functions $\varphi',\varphi'',\psi',\psi''\in
C(Y)$ such that

\item $\varphi'f\le\varphi_2\le \varphi''f$;

\item$\psi'f\le\psi_2\le \psi''f$;
\item $\varphi'|B=\varphi''|B=\psi|B=\psi''|B=\min\{\inf\varphi,\inf\psi\}-1$.
\end{enumerate}

Then we have $$I(f)(\mu)(\varphi')=\mu(\varphi'f)\le
\mu(\varphi_2)\le \mu(\varphi''f)I(f)(\mu)(\varphi'')$$ and, since
$$I(f)(\mu)(\varphi')=I(f)(\mu)(\varphi'')=\min\{\inf\varphi,\inf\psi\}-1,$$
we conclude that $\mu(\varphi_2)=\min\{\inf\varphi,\inf\psi\}-1$.
Similarly, one can show that
$\mu(\psi_2)=\min\{\inf\varphi,\inf\psi\}-1$. We have obtained a
contradiction.
\end{proof}

Since the functor $I$ preserves preimages,  for any $f\colon X\to
Y$ and $\mu\in I(X)$, we have
$\mathrm{supp}(I(f)(\mu))=f(\mathrm{supp}(\mu))$. This follows
from general properties of functors in the category $\comp$
established in \cite{S}.

A functor $F$ in the category $\comp$ is called {\em normal} (see
\cite{S}) if $F$ is continuous, preserves weight, singletons,
empty set, the onto maps, embeddings, intersections, and
preimages.

\begin{prop} The functor $I$ is normal.
\end{prop}
\begin{proof} This follows from Propositions \ref{p:emb}, \ref{p:int}, \ref{p:onto}, \ref{p:prei}, and also from
the fact that $I$ is a subfunctor of $O$ and the latter functor is
almost normal in the sense that it satisfies all the properties
from the definition of the normal functor excepting the
preimage-preserving property (see \cite{R}). As an example, we
remark that the functor $O$ preserves the weight of infinite
compact Hausdorff spaces, i.e. $w(X)=w(O(X))$, for any infinite
$X$. Since $I(X)\subset O(X)$, we conclude that $w(X)=w(I(X))$.

\end{proof}

\begin{prop}\label{p:fins} Let $|X|=n$, then the space $I(X)$ is homeomorphic to the
$(n-1)$-dimensional simplex.
\end{prop}
\begin{proof} Let $X=\{x_1,\dots,x_n\}$. We first show that, for every $\mu\in I(X)$, there
exist $\lambda_1,\dots,\lambda_n\in\R$ such that
$\mu(\varphi)=\max\{\varphi(x_i)+\lambda_i\mid i=1,\dots,n\}$, for
every $\varphi\in C(X)$. For every $i=1,\dots,n$, define
$\varphi_i$ by the formula $\varphi_i(x_j)=\delta_{ij}-1$. Let
$\lambda_i=\inf\{\mu(\alpha\varphi_j)\mid \alpha\ge0\}$.

Now, let $\varphi\in C(X)$. Then, for every $\alpha\ge0$, we have
$$\varphi\le \max\{\alpha\varphi_i+\varphi(x_i)\mid i=1,\dots,n\},$$
whence
$$\mu(\varphi)\le\max\{\varphi(\alpha\varphi_i)+\varphi(x_i)\mid
i=1,\dots,n\}.$$ Passing to the limit as $\alpha\to\infty$, we see
that $\mu(\varphi)\le\max\{\lambda_i+\varphi(x_i)\mid
i=1,\dots,n\}$. Actually, for sufficiently large $\alpha$, we have
the equality. The measure $\mu$ defined above is denoted by
$\oplus_{i=1}^n\lambda_i\odot\delta_{x_i}$.

Let $\Gamma^{n-1}=\{(\lambda_1,\dots,\lambda_n)\in\R^n\mid
\max\{\lambda_1,\dots,\lambda_n\}=0\}$. It is evident that
$\Gamma^{n-1}$ is homeomorphic to the $(n-1)$-dimensional simplex.
Define the map $\xi\colon \Gamma^{n-1}\to I(\{x_1,\dots,x_n\})$ by
the formula $\xi(\lambda_1,\dots,\lambda_n)=
\oplus_{i=1}^n\lambda_i\odot\delta_{x_i}$. From what was already
prove it follows that  $\xi$ is an onto map.

Given distinct
$(\lambda_1,\dots,\lambda_n),(\lambda'_1,\dots,\lambda'_n)\in\Gamma^{n-1}$,
one can find $i$ such that $\lambda_i\neq\lambda'_i$. Define
$\varphi\colon X\to\mathbb R$ as follows:
$$\varphi(x)=\begin{cases} -1-\max\{\lambda_j,\lambda'_j\},& \text{ if }x=x_j,\ i\neq
j,\\
-\max\{\lambda_i,\lambda'_i\}& \text{ if }x=x_i.
\end{cases}
$$
Then
$$\xi(\lambda_1,\dots,\lambda_n)=\lambda_i-\max\{\lambda_i,\lambda'_i\}\neq
\lambda'_i-\max\{\lambda_i,\lambda'_i\}=\xi(\lambda'_1,\dots,\lambda'_n),$$
whence we conclude that the map $\xi$ is also an embedding.

Since the set  $\Gamma^{n-1}$ is compact, in order to show that
$\xi$ is a homeomorphism, one has to verify that the map $\xi$ is
continuous.

Let $\mu=\xi(\lambda_1,\dots,\lambda_n)$ and
$\langle\mu;\varphi;\varepsilon\rangle$ be a subbase neighborhood
of $\mu$ in $I(X)$, where $\varphi\in C(X)$, $\varepsilon>0$.
Define neighborhoods $U_i$ of $\lambda_i$, $i=1,\dots,n$, in $\R$
as follows:
$$U_i=\begin{cases} (\lambda_i-\varepsilon,\lambda_i+\varepsilon,
& \text{if }\lambda_i>-\infty,\\
[-\infty,\min\{\lambda_j\mid \lambda_j>-\infty,\
j=1,\dots,n\}-\varepsilon),& \text{if }\lambda_i=-\infty.
\end{cases}$$
Then $U=(U_1\times\dots\times U_n)\cap \Gamma^n$ is a neighborhood
of $(\lambda_1,\dots,\lambda_n)$ with $\xi(U)\subset
\langle\mu;\varphi;\varepsilon\rangle$ and thus the map $\xi$ is
continuous.

\end{proof}

\begin{prop} The set
$$I_\omega(X)=\{\oplus_{i=1}^n\lambda_i\odot\delta_{x_i}\mid
\lambda_i\in\R,\ i=1,\dots,n,\ \oplus_{i=1}^n\lambda_i=0,\ x_i\in
X,\  n\in\mathbb N \}$$ (i.e., the set of idempotent probability
measures of finite support) is dense in $I(X)$.
\end{prop}
\begin{proof} It follows from Proposition  \ref{p:onto} and
results of general theory of functors in the category $\comp$ (see
\cite{S}) that the set $I_\omega(X)$ of the idempotent measures
with finite supports is dense in $I(X)$. The  statement is now a
consequence of Proposition \ref{p:fins}.
\end{proof}

We see that the spaces $I(X)$ and $P(X)$ are homeomorphic for
every finite $X$. In forthcoming publications we will show that
they are also homeomorphic for infinite metrizable $X$. However,
the following statement holds.

\begin{prop} The functors $P$ and $I$ are not isomorphic.
\end{prop}
\begin{proof} Let $X=\{a,b,c\}$, $Y=\{a,b\}$, $Z=\{a,c\}$, where
$a ,b,c$ are distinct points. Denote by $f\colon X\to Y$ and
$g\colon X\to Z$ the retractions such that $f(c)=b$ and $g(b)=c$.
Then the map $(P(f),P(g))\colon P(X)\to P(Y)\times P(Z)$ is
obviously an embedding while the map $(I(f),I(g))\colon I(X)\to
I(Y)\times I(Z)$ is not. Indeed, let
\begin{align*}\mu=&(-1)\odot\delta_a\oplus0\odot\delta_b\oplus0\odot\delta_c,\\
\nu=&(-2)\odot\delta_a\oplus0\odot\delta_b\oplus0\odot\delta_c,\end{align*}
then $$I(f)(\mu)=I(f)(\nu)=0\odot\delta_a\oplus0\odot\delta_b,\
I(g)(\mu)=I(g)(\nu)=0\odot\delta_a\oplus0\odot\delta_c.$$

\end{proof}

Given $x,y\in\mathbb R^n$ and $\lambda\in\mathbb R$, we denote by
$x\oplus y$ the coordinatewise maximum of $x$ and $y$ and by
$\lambda\odot x$ the vector obtained from $x$ by adding $\lambda$
to every its coordinate. A subset $A$ in $\mathbb R^n$ is called
{\em max-plus convex} if $\lambda_1\odot x_1\oplus\lambda_2\odot
x_2\in A$ whenever $x_1,x_2\in A$ and $\lambda_1,\lambda_2\in
\mathbb R$ with $\lambda_1\oplus\lambda_2=0$. Note that, in this
definition, one can also assume that $\lambda_1,\lambda_2\in \R$.

One can similarly define, for any $\mu_1,\mu_2\in I(X)$ and
$\lambda_1,\lambda_2\in \R$ with $\lambda_1\oplus\lambda_2=0$, the
{\em max-plus convex combination} $\mu=\lambda_1\odot
\mu_1\oplus\lambda_2\odot \mu_2$ as follows:
$$\mu(\varphi)=\lambda_1\odot \mu_1(\varphi)\oplus\lambda_2\odot
\mu_2(\varphi),\ \varphi\in C(X).$$

The following statement is obvious.
\begin{prop} We have $\lambda_1\odot
\mu_1\oplus\lambda_2\odot \mu_2\in I(X)$.
\end{prop}

\begin{prop} Let $A\subset I(X)$, $A\neq\emptyset$. Then $\sup
A\in I(X)$.
\end{prop}

Let $A\subset \mathbb R^n$ be a compact max-plus convex subset. By
abusing the language, we denote by $x_1,\dots,x_n$ the coordinate
functions $\mathbb R^n\to\mathbb R$. Given $\mu\in I(A)$, we let
$\beta_A(\mu)=(\mu(x_1),\dots,\mu(x_n))$.

\begin{prop}\label{p:bary} The map $\beta=\beta_A\colon I(A)\to A$ is
continuous.
\end{prop}
\begin{proof} The continuity of the map $\mu\mapsto\beta_X(\mu)$  follows from the fact
that if $\mu'\in \langle\mu;x_1,\dots,x_n;\varepsilon\rangle$,
then $\|\beta(\mu)-\beta(\mu')\|<\varepsilon$.

Given $\mu=\oplus_{i=1}^k\lambda_i\odot\delta_{a_i}\in I(A)$, we
see that $\beta(\mu)= \oplus_{i=1}^k\lambda_i\odot a_i\in A$.
Since $I_\omega(X)$ is dense in $I(X)$, we see that $\beta(\mu)\in
A$ for every $\mu\in I(A)$. Therefore, the map $\beta$ is
well-defined.

\end{proof}

The map $\beta\colon I(A)\to A$ is called the {\em idempotent
barycenter map}.
\begin{rem} One can extend Proposition \ref{p:bary} over the case
of the compact max-convex subsets in arbitrary Tychonov power
$\mathbb R^\tau$.
\end{rem}

\section{Idempotent probability measure monad}

It is known \cite{R} that the functor of order-preserving
functionals forms a monad on the category $\comp$. In this section
we are going to show that the functor $I$ is the functorial part
of a submonad of the monad of order-preserving functionals.

A {\it monad} $\mathbb T=(T,\eta,\mu)$ in the category ${\mathcal
E}$ consists of an endofunctor $T\colon {\mathcal E}\to{\mathcal
E}$ and natural transformations $\eta\colon 1_{\mathcal E}\to T$
(unity), $\mu\colon T^2\to T$ (multiplication) satisfying the
relations $\mu\circ T\eta=\mu\circ\eta T=${\bf 1}$_T$ and
$\mu\circ\mu T=\mu\circ T\mu$.

A natural transformation $\psi\colon T\to T'$ is called a {\it
morphism} from a monad $\mathbb T=(T,\eta,\mu)$ into a monad
$\mathbb T'=(T',\eta',\mu')$ if $\psi\circ\eta= \eta'$ and
$\psi\circ\mu=\mu'\circ\eta T'\circ T\psi$. If all the components
of $\psi$ are monomorphisms then the monad $\mathbb T$ is called a
{\it submonad} of $\mathbb T'$.

If  $\mathbb T=(T,\eta,\mu)$ is a monad in the category ${\mathcal
E}$, then a pair $(X,\xi)$, where $\xi\colon T(X)\to X$ is a {\em
$\mathbb T$-algebra} if $\xi\eta_X=\id_X$ and $\xi\mu_X=\xi
T(\xi)$. Given $\mathbb T$-algebras $(X,\xi)$, $(X',\xi')$, we say
that a morphism $f\colon X\to X'$ is a {\em morphism of $\mathbb
T$-algebras} if $f\xi=\xi'T(f)$. The $\mathbb T$-algebras and
their morphisms form a category.

The {\em hyperspace monad } $\mathbb H=(\exp,s,u)$ in the category
$\comp$ is defined as follows. The natural transformation $s\colon
\id\to \exp$ acts by the formula $s_X(x)=\{x\}$, $x\in X$. The
natural transformation $u\colon\exp^2\to\exp$ is defined by the
formula $u(\mathcal A)=\cup\mathcal A$, $\mathcal A\in\exp^2 X$.

Let $X\in|\comp|$. Given $\varphi\in C(X)$, define
$\bar\varphi\colon I(X)\to\mathbb R$ as follows:
$\bar\varphi(\mu)=\mu(\varphi)$, $\mu\in I(X)$.

\begin{lemma}\label{l:1} If $\varphi\in C(X)$ and $\lambda\in\R$, then
$\overline{\lambda\odot\varphi}=\lambda\odot\bar\varphi$.
\end{lemma}
\begin{proof} Given $\mu\in I(X)$, we have $\overline{\lambda\odot\varphi}(\mu)=
\mu(\lambda\odot\varphi)=\lambda\odot\mu(\varphi)=\lambda\odot\bar\varphi(\mu)$.
\end{proof}

\begin{lemma}\label{l:2}  If $\varphi,\psi\in C(X)$, then
$\overline{\varphi\oplus\psi}=\bar\varphi\oplus\bar\psi$.
\end{lemma}
\begin{proof} Given $\mu\in I(X)$, we have $\overline{\varphi\oplus\psi}(\mu)=\mu(\varphi\oplus\psi)=
\mu(\varphi)\oplus\mu(\psi)=\bar\varphi(\mu)\oplus\bar\psi(\mu)=(\bar\varphi\oplus\bar\psi)(\mu)$.
\end{proof}

Given $M\in I^2(X)$, define the map $\zeta_X(M)\colon
C(X)\to\mathbb R$ as follows:
$\zeta_X(M)(\varphi)=M(\bar\varphi)$.

\begin{prop}  We have $\zeta_X(M)\in I(X)$.
\end{prop}
\begin{proof}
Check the conditions from the definition of $I(X)$.

 1.
$\zeta_X(M)(c_X)=M(\overline{c_X})= M(c_{I(X)})=c$.

2. Applying Lemma \ref{l:1} we obtain
$\zeta_X(M)(\lambda\odot\varphi)=M(\overline{\lambda\odot\varphi})=M(\lambda\odot\bar\varphi)=\lambda\odot
M(\bar\varphi)=\lambda\odot \zeta_X(M)(\varphi)$.

3. Applying Lemma \ref{l:2} we obtain
$\zeta_X(M)(\varphi\oplus\psi)=M(\overline{\varphi\oplus\psi})=M(\bar\varphi\oplus\bar\psi)=M(\bar\varphi)\oplus
M(\bar\psi)=\zeta_X(M)(\varphi)\oplus\zeta_X(M)(\psi)$.

\end{proof}

Thus, we obtain a map $\zeta_X\colon I^2(X)\to I(X)$. It follows
from the results of \cite{R} that $\zeta=(\zeta_X)$ is a natural
transformation from the functor $I^2$ to the functor $I$.
Actually, this natural transformation is the restriction of the
natural transformation $O^2\to O$ defined by Radul.

\begin{thm} The triple $\mathbb I=(I,\eta,\zeta)$ is a monad on the
category $\comp$.
\end{thm}
\begin{proof}
As we already remarked, the natural transformation $\zeta$ is the
restriction of the natural transformation $O^2\to O$ and $\delta$
maps the identity functor into $I\subset O$. Therefore, the result
follows from the fact that the functor $O$ generates a monad in
$\comp$ (see \cite[Theorem 3]{R}).

\end{proof}

Actually, $\mathbb I$ is a submonad of the monad $\mathbb O$
generated by the functor $O$ (see \cite{R}).

\begin{prop}\label{p:alg} Let $X$ be a  compact max-plus convex subset in
$\mathbb R^n$. Then the pair $(X,\beta)$ is an $\mathbb
I$-algebra.
\end{prop}
\begin{proof} Let $x\in X$, then
$\beta\eta_X(x)=\beta(\delta_x)=x$, i.e. $\beta\eta_X=\id_X$.

Let $M=\oplus_{i=1}^n\lambda_i\odot\mu_i\in I^2(X)$, where
$\mu_i=\oplus_{j=1}^m\kappa_{ij}\odot\delta_{x_{ij}}$. Then
\begin{align*}\beta\zeta_X(M)=&\beta(\oplus_{i=1}^n\lambda_i\odot\mu_i)=
\oplus_{i=1}^n\lambda_i\odot\left(\oplus_{j=1}^m\kappa_{ij}\odot\delta_{x_{ij}}\right)
\\
=&\oplus_{i=1}^n\oplus_{j=1}^m\left(\lambda_i\odot\kappa_{ij}\right)\odot\delta_{x_{ij}}
\\ = &\beta (I(\beta))(M).\end{align*}

Since the idempotent probability measures $M$ of the form as above
are dense in $I^2(X)$, we are done.
\end{proof}

A continuous map $f\colon A\to B$ of compact max-plus convex sets
in euclidean space is called {\em max-plus affine} if it preserves
the max-plus convex combinations (i.e. $f(\alpha\odot
a\oplus\beta\odot b)= \alpha\odot f(a)\oplus \beta\odot f(b)$, for
every $a,b\in A$ and $\alpha,\beta\in\R$ with
$\alpha\oplus\beta=0$). Examples of max-plus affine maps are the
projections onto coordinate hyperplanes in $\mathbb R^n$.

\begin{prop} Let $f\colon X\to Y$ be a continuous affine map of
max-plus compact convex subsets in euclidean spaces. Then $f$ is a
morphism of $\mathbb I$-algebras.
\end{prop}

\begin{proof} One has to show that
$f\beta_X(\mu)=\beta_YI(f)(\mu)$, for any $\mu\in I(X)$. It
suffices to verify this equality for $\mu$ of finite support. If
$\mu=\oplus_{i=1}^n\lambda_i\odot\delta_{x_i}$, then
$$f\beta_X(\mu)=f(\oplus_{i=1}^n\lambda_i\odot
x_i)=\oplus_{i=1}^n\lambda_i\odot
f(x_i)=\beta_Y(\oplus_{i=1}^n\lambda_i\odot\delta_{f(x_i)})=\beta_YI(f)(\mu).$$
\end{proof}

\begin{prop} The diagonal map $$\Phi=(\bar\varphi)_{\varphi\in
C(X)}\colon I(X)\to \prod\{\mathbb R_\varphi=\mathbb R\mid
\varphi\in C(X)\}=\mathbb R^{C(X)}$$ embeds $I(X)$ as a max-plus
convex subset of $\mathbb R^{C(X)}$.
\end{prop}
\begin{proof} This follows from the equality
$$\bar\varphi(\alpha_1\odot\mu_1\oplus\alpha_2\odot\mu_2) =
\alpha_1\odot\bar\varphi(\mu_1)\oplus\alpha_2\odot\bar\varphi(\mu_2),$$
for any $\mu_1,\mu_2\in I(X)$, $\varphi\in C(X)$, and
$\alpha_1,\alpha_2\in\R$ with $\alpha_1\oplus\alpha_2=0$.
\end{proof}

Using the monad structure for the functor $I$, one can define, for
all $\mu\in I(X)$, $\nu\in I(Y)$, the {\em tensor product}
$\mu\otimes\nu\in I(X\times Y)$ (see, e.g. \cite{Za}). For the
sake of completeness, we recall its construction. For every $y\in
Y$, let $i_y\colon X\to X\times Y$ denote the map defined by the
formula $i_y(x)=(x,y)$, $x\in X$. Then define the map
$g_{\mu}\colon Y\to I(X\times Y)$ by the formula
$g_{\mu}(y)=I(i_y)(\mu)$, $y\in Y$. Finally,
$$\mu\otimes\nu=\zeta_{X\times Y}I(g_{\mu})(\nu).$$

If $\mu=\oplus_{i=1}^m\lambda_i\odot\delta_{x_i}\in I(X)$,
$\nu=\oplus_{j=1}^n\kappa_j\odot\delta_{y_j}\in I(Y)$, then
$$\mu\otimes\nu=\bigoplus_{i=1}^m \bigoplus_{j=1}^n (\lambda_i\odot\kappa_j)\odot \delta_{(x_i,y_j)}\in I(X\times Y).$$

By induction, one can define the tensor product for arbitrary
finite products: if $\mu_i\in I(X_i)$, $i=1,\dots,n$, then
$$\mu_1\otimes\dots\otimes\mu_n=(\mu_1\otimes\dots\otimes\mu_{n-1})\otimes\mu_n\in
I((X_1\times\dots\times X_{n-1})\times X_n)=I(X_1\times\dots\times
X_{n}).$$

The tensor product can be also defined for infinite products of
idempotent probability measures. Given $\mu_\alpha\in
I(X_\alpha)$, $\alpha\in A$, where $A$ is an arbitrary infinite
set of indices, one defines $\otimes\{\mu_\alpha\mid \alpha\in
A\}$ as a unique $\nu\in I(\prod_{\alpha\in A}X_\alpha)$
satisfying the property $I(\pr_B)(\nu)=\otimes\{\mu_\alpha\mid
\alpha\in B\}$, for every nonempty finite subset $B$ of $A$.

Note that the probability measure monad $\mathbb P$ is also a
submonad in $\mathbb O$. The intersection of the submonads
$\mathbb P$ and $\mathbb I$ (in a natural sense) is the identity
submonad of the monad $\mathbb O$.

\begin{thm} The hyperspace monad $\mathbb H$ is a submonad of the
monad $\mathbb I$.
\end{thm}
\begin{proof} Given a compact Hausdorff space $X$, define a map
$j_X\colon \exp X\to I(X)$ by the condition:
$j_X(A)(\varphi)=\max(\varphi|A)$. It is straightforward to verify
that $j_X$ is well-defined.

We are going to demonstrate that the map $j_X$ is continuous. Let
$A_0\in\exp X$ and $\langle j_X(A_0);\varphi;\varepsilon\rangle$
be a subbase neighborhood of $j_X(A_0)$ in $I(X)$. There exists a
finite open in $X$ cover $\mathcal U=\{U_1,\dots,U_n\}$ of $A_0$
such that, for every $U_i\in\mathcal U$, the oscillation of
$\varphi$ on $U_i$ (i.e. the number
$|\sup(\varphi|U_i)-\inf(\varphi|U_i)|$) is less than
$\varepsilon$. Let $A\in\langle U_1,\dots,U_n\rangle$. We have to
show that $j_X(A)\in \langle j_X(A_0);\varphi;\varepsilon\rangle$.

Let $j_X(A)(\varphi)=\varphi(a)$, where $a\in U_i\cap A$, for some
$i$. Then there is $a_0\in U_i\cap A_0$ and
$|\varphi(a)-\varphi(a_0)|<\varepsilon$, whence
$j_X(A)(\varphi)=\varphi(a)<\varphi(a_0)+\varepsilon\le
j_X(A_0)(\varphi)+\varepsilon$. Proceeding similarly, we prove
that $j_X(A_0)(\varphi)< j_X(A)(\varphi)+\varepsilon$.

Note that the map $j_X$ is an embedding. Indeed, let $A,B\in \exp
X$ and $A\neq B$. Without loss of generality, we may assume that
$A\setminus B\neq\emptyset$. Let $\varphi \in C(X)$ be a function
with  the following properties: $\varphi|B\equiv 0$,
$\varphi(x)>0$, for some $x\in A\setminus B$. Then
$j_X(A)(\varphi)>j_X(B)(\varphi)=0$.

Given $f\colon X\to Y$ and $A\in\exp X$, $\varphi\in C(Y)$, we see
that
\begin{align*} (I(f)j_X(A))(\varphi)=&  j_X(A)(\varphi f)=\max\{\varphi f(a)\mid a \in
A\}\\ =&\max\{\varphi (b)\mid b \in f(A)\}=j_Y(f(A))(\varphi),
\end{align*}
whence $I(f)j_X=j_Y\exp f$ and we see that $j\colon \exp\to I$ is
a natural transformation.

We are going to prove that $j$ is a monad morphism. To this end,
show that the
diagram $$\xymatrix{\exp^2X\ar[rr]^{I(j_X)j_{\exp X}}\ar[d]_{u_X}&&I^2(X)\ar[d]^{\zeta_X}\\
\exp X\ar[rr]^{j_X}&& I(X)}$$ is commutative. We  prove this for
points of with finite supports. Let $\mathcal A\in\exp^2X$,
$\mathcal A=\{A_1,\dots,A_k\}$, where
$A_i=\{a_{i1},\dots,a_{il}\}$.

Then $j_{\exp X}(\mathcal A)=\oplus_{p=1}^k 0\odot\delta(A_p)$ and
$$ I(j_X)j_{\exp X}(\mathcal A)=
I(j_X)(\oplus_{p=1}^k0\odot\delta(A_p))=\oplus_{p=1}^k0\odot\delta(\oplus_{q=1}^l0\odot\delta(a_{pq})),
$$
and
$$ \zeta_X I(j_X) j_{\exp X}(\mathcal A)=\oplus_{p=1}^k\oplus_{q=1}^l(0\odot0)\odot\delta(a_{pq}).$$

On the other hand, $$j_Xu_X(\mathcal A)=j_X(\{a_{pq}\mid 1\le p\le
k,\ 1\le q\le
l\})=\oplus_{p=1}^k\oplus_{q=1}^l0\odot\delta(a_{pq}).
$$
Since the points of finite support are dense in $\exp^2X$, we are
done.

Also $j_Xs_X(x)=j_X(\{x\})=\delta_x$, for every $x\in X$, and we
see that $J$ is a monad morphism.
\end{proof}
\begin{rem} Let $\chi\colon X\to [0,1]$ be a  fuzzy set such that the map $\chi$ is continuous
and  $\chi^{-1}(1)\neq\emptyset$.  One can identify $\chi$ with an
element $j_X(\chi)$ of $I(X)$ as follows:
$$j_X(\chi)(\varphi)=\sup\{\varphi(x)+\ln\chi(x)\mid x\in X\},\ \varphi\in C(X).$$
\end{rem}

\section{Milyutin maps of idempotent probability measures}

The notion of Milyutin map was first introduced for the
probability measure functor (see, e.g., \cite{P} for the
construction).

\begin{thm}\label{t:mil} Let $X$ be a compact metrizable space. Then there
exists a zero-dimensional compact metrizable space $X$ and a
continuous map $f\colon X\to Y$ for which there exists a
continuous map $s\colon Y\to I(X)$ such that
$\mathrm{supp}(y)\subset f^{-1}(y)$, for every $y\in Y$.
\end{thm}
\begin{proof} One can easily construct a sequence $(\mathcal
W_i)$, where each $\mathcal W_i$ is a finite set of pairs of
subsets of $Y$ satisfying the properties:
\begin{enumerate}
\item $\mathcal U_i=\{U\mid (U,V)\in\mathcal W_i\}$ and $\mathcal V_i=\{U\mid (U,V)\in\mathcal W_i\}$
are finite closed covers of the space $Y$;
\item $U\subset\mathrm{Int}_Y(V)$ for every $(U,V)\in \mathcal W_i$;
\item $\mathrm{mesh}(V_i)<(1/i)$ for every $i$ (we assume that some metric is fixed on $Y$; the mesh of a
family of subsets in a metric space is the supremum of the
diameters of its members).
\end{enumerate}

We let $X_i=\coprod\{V\mid (U,V)\in \mathcal W_i\}$. The map
$f_i\colon X_i\to Y$ is the map such that $f_i|V\colon V\to Y$ is
the inclusion map for every $V$ such that $(U,V)\in \mathcal W_i$.
Let $\alpha_i\colon X_i\to[-\infty,0]$ be a continuous function
such that, for every $(U,V)\in \mathcal W_i$, we have
$\alpha_i|U\equiv0$ and $\alpha_i|(V\setminus
\mathrm{Int}_Y(V))\equiv-\infty$.

Let $$X=\left\{(x_i)_{i=1}^\infty \in \prod_{i=1}^\infty X_i\mid
f_i(x_i)=f_j(x_j)\text{ for every }i,j \right\}.$$ Define the map
$f\colon X\to Y$ by the formula $f((x_i)_{i=1}^\infty )=f_1(x_1)$.

Given $y\in Y$, define
$$s(y)=\bigotimes_{i=1}^\infty\oplus\{\alpha_i(x)\odot\delta_x\mid
x\in f_i^{-1}(y)\}\in I(X).$$ It is easy to see that $s$ is
well-defined and continuous. For any $y\in Y$, we have
$\mathrm{supp}(s(y))=\prod_{i=1}^\infty f^{-1}(y)$ and therefore
$f(s(y))=\delta_y$, for every $y\in Y$.

Finally, we leave to the reader the verification that the space
$X$ is zero-dimensional and compact metrizable.

\end{proof}

Similarly as in the case of probability measures, one can show
that the product of idempotent Milyutin maps is Milyutin and than
the restriction of a Milyutin map onto a full preimage of a closed
set is also Milyutin. This allows us to prove that every compact
Hausdorff space is the image of a zero-dimensional compact
Hausdorff space under a Milyutin map.

We call a map $f\colon X\to Y$ that satisfies the properties of
Theorem \ref{t:mil} an {\em idempotent Milyutin map}.

We will need the following notion introduced by E. Shchepin
\cite{Sh}. A commutative diagram
\begin{equation}\label{dia:1} \xymatrix{X\ar[r]^f\ar[d]_g &Y\ar[d]^u\\
Z\ar[r]^v & T}
\end{equation} is called {\it bicommutative}
 if its
{\it characteristic map} $$\chi=(f,g)\colon X\to
Y\times_TZ=\{(y,z)\in Y\times Z\mid u(y)=v(z)\}$$ is onto.

\begin{thm}\label{t:open} The idempotent probability measure functor is open.
\end{thm}
\begin{proof} We first consider the case of surjective map of
finite spaces. Let $f\colon X\to Y$ be such a map. Since the
composition of any two open maps is open, without loss of
generality, one may assume that $X=\{x_0,x_1,\dots,x_n\}$,
$Y=\{y_1,\dots,y_n\}$, and the map $f\colon X\to Y$  acts by the
formula $f(x_0)=y_1$, $f(x_i)=y_i$, $y=1,\dots,n$. Let $\mu_0\in
I(X)$, $\mu_0=\oplus_{i=0}^n\alpha_{i0}\odot\delta_{x_i}$,
$\nu_0=I(f)(\mu_0)$, and $(\nu_k)_{k=1}^\infty$ be a sequence in
$I(Y)$ converging to $\nu_0$. We have
$\nu_k=\oplus_{j=1}^n\beta_{kj}\odot\delta_{y_j}$. Then
$\lim_{k\to\infty}\beta_{kj}=\alpha_{0j}$, for $j=2,\dots,n$, and
$\lim_{k\to\infty}\beta_{k1}=\max\{\alpha_{00},\alpha_{01}\}$.
Without loss of generality, we may assume that
$\alpha_{00}\ge\alpha_{01}\}$. Then let $\alpha_{k0}=\beta_{k1}$,
$\alpha_{k1}=\min\{\beta_{k1},\alpha_{01}\}$. Let
$\mu_k=\oplus_{i=0}^n\alpha_{ik}\odot\delta_{x_i}$, $k\in\mathbb
N$. It is obvious that $I(f)(\mu_k)=\nu_k$, for every $k$, and
$\lim_{k\to\infty}(\mu_k)=\mu_0$. This is equivalent to the
openness of the map $I(f)$.

Let $C$ denote the Cantor set. Let us prove that the map $I(\pr)$,
where $\pr\colon C\times C\to C$ denotes the projection onto the
first factor, is open. To this end, represent $C$ as
$\varprojlim\{C_i,f_{ij}\}$, where $C_i$ are finite sets and
$f_{ij}\colon C_i\to C_j$ are surjections, $i\ge j$. From the
results of \cite{S} it follows that, in order to prove that
$I(\pr)$ is open, it is sufficient to prove that the diagram
$$\xymatrix{I(C_i\times C_i)\ar[d]_{I(\pi_i)}\ar[rr]^{I(f_{ij}\times f_{ij})}& &I(C_j\times C_j)\ar[d]^{I(\pi_j)}\\
I(C_i)\ar[rr]^{I(f_{ij})}& &I(C_j)}$$ (here $\pi_k\colon C_k\times
C_k\to C_k$ denotes the projection onto the first factor) is
bicommutative i.e. the map
\begin{align*} (I(\pi_i),I(f_{ij}\times f_{ij}))\colon &I(C_i\times C_i)\to
I(C_j\times C_j)\times_{I(C_j)} I(C_i) \\ =&\{(\mu,\nu)\in
I(C_i)\times I(C_j\times C_j)\mid I(f_{ij})(\mu)=I(\pi_j)(\nu)\}
\end{align*} (called
the characteristic map of the diagram) is an onto map.

Without loss of generality, one may assume that
$$C_j=\{x_1,\dots,x_p\},\ C_i=\{y_0,y_1,\dots,y_p\}$$ (all the points are assumed to be distinct) and the map
$f_{ij}$ act as follows: $f_{ij}(y_m)=x_m$, $m=1,\dots,p$,
$f_{ij}(y_0)=y_1$.  Thus, given $(\mu,\nu)\in I(C_j\times
C_j)\times_{I(C_j)} I(C_i)$, one can write
$$\mu=\bigoplus_{k=0}^p\kappa_k\odot\delta_{y_k},\
\nu=\bigoplus_{m,n=1}^p\lambda_{mn}\odot\delta_{(x_m,x_n)}.$$
Without loss of generality, we may assume that
$\kappa_0\le\kappa_1$. Define
$$\nu'=\bigoplus_{m,n=0}^p\lambda'_{mn}\odot\delta_{(y_m,y_n)}\in
I(C_i\times C_i)$$ by the conditions $\lambda'_{mn}=\lambda_{mn}$,
for $m\ge1$, $n=0,1,\dots,p$,
$\lambda'_{0n}=\min\{\kappa_0,\lambda_{1n}\}$, $n=0,1,\dots,p$. We
leave to the reader the verification of the fact that $\nu'$ is as
required.

Now, consider an open map $f\colon X\to Y$ of compact metrizable
spaces. Let $p\colon Z\to Y$ be an idempotent Milyutin map, where
$Z$ is a compact metrizable zero-dimensional space. We may assume
that $Z$ is homeomorphic to the Cantor set, $Z$ is a subset of the
product $T\times Y$ and $p$ coincides with the restriction of the
projection $\hat p\colon T\times Y\to Y$ onto the second factor.

Denote by $p\colon Z\times_YX\to X$ the projection map,
$p(z,x)=x$. We assume that $Z\times_YX\subset T\times Y\times X$
For every $x\in X$, let $i_x\colon T\times Y \to T\times Y\times
X$ be the map defined by the formula $i_x(t,y)=(t,y,x)$.

Let $s\colon Y\to I(Z)$ be a map such that $\hat p s(y)=\delta_y$,
for every $y\in Y$. Now, consider a sequence $(\nu_i)$ in $I(Y)$
converging to $\nu_0$ and $\mu_0\in I(X)$ such that
$I(f)(\mu_0)=\nu_0$. Define a map $g\colon X\to I(T\times Y\times
X)$ by the formula $g(x)=I(i_x)(s(f(x)))$, $x\in X$. Let $\mu_0'=
\zeta_{T\times Y\times X}(I(g)(\mu_0))$.

Denote by $\pi_i$ the projection of $T\times Y\times X$ onto the
$i$-th factor and by $\pi_{ij}$ the projection of $T\times Y\times
X$ onto the product of the $i$-th and $j$-th factors. We then have
$$I(\pi_3)(\mu_0')=I(\pi_3)\zeta_{T\times Y\times
X}(I(g)(\mu_0))=\zeta_XI^2(\pi_3)(I(g)(\mu_0))=\zeta_XI(\eta_X)(\mu_0)=\mu_0.$$

For every $i=0,1,2,\dots$, define $\nu'_i=\zeta_Z I(s)(\nu_i)$.
Then $$I(p)(\nu'_i)=I(p)\zeta_Z I(s)(\nu_i)=\zeta_Y
I^2(p)I(s)(\nu_i)=\zeta_Y I(I(p)s)(\nu_i)= \zeta_Y
I(\eta_Y)(\nu_i)=\nu_i.$$

We have
\begin{align*} I(\pi_{12})(\mu'_0)=&I(\pi_{12})\zeta_{T\times
Y\times X}(I(g)(\mu_0))=\zeta_{T\times
Y}I^2(\pi_{12})(I(g)(\mu_0))\\=&\zeta_{T\times
Y}I(I(\pi_{12}g))(\mu_0)=\zeta_{T\times
Y}I(sf)(\mu_0)=\zeta_{T\times Y}I(s)I(f)(\mu_0)\\ =&
\zeta_{T\times Y}I(s)(\nu_0)=\nu'_0.
\end{align*}

Let $h\colon K\to \pi_{12}(Z)$ be an open onto map of a
zero-dimensional compact metrizable space. Without loss of
generality, one may assume that the composition $\pi_{12}h$ is
homeomorphic to the projection map  $\pr\colon C\times C\to C$.
Let $\mu_0''\in I(K)$ be such that $I(h)(\mu_0'')=(\mu_0')$. Then,
by the openness of the map $I(\pi_{12}h)$, there exists a sequence
$(\mu''_i)$ in $I(K)$ such that $\lim_{i\to\infty}\mu''_i=\mu_0''$
and $I(\pi_{12}h)(\mu''_i)=\nu'_i$.

Let $\mu_i=I(\pi_3h)(\mu_i'')$. Then
$$\lim_{i\to\infty}\mu_i=\lim_{i\to\infty}I(\pi_3h)(\mu''_i)=I(\pi_3h)(\mu_0'')=\mu_0.$$
For every $i\in\mathbb N$, we have
$$I(f)(\mu_i)=I(f\pi_3h)(\mu_i'')=I(\hat
p\pi_{12}h)(\mu_i'')=I(\hat p)(\nu'_i)=\nu_i.$$ This proves that
$I(f)$ is an open map.

\end{proof}

A functor $F$ in the category $\comp$ is called {\em
bicommutative} if $F$ preserves the class of bicommutative
diagrams.

\begin{cor} The functor $I$ is bicommutative.
\end{cor}
\begin{proof} The fact follows from Theorem \ref{t:open} and the
result due to Shchepin \cite{S} that every open functor is
bicommutative.
\end{proof}

\section{Correspondences of idempotent probability measures with restricted marginals}

Given a finite collection $X_1,\dots, X_k$ of compact Hausdorff
spaces, define a map $M_{X_1,\dots,X_k}\colon I(\prod X_i)\to
\prod I(X_i)$ as follows:
$$M_{X_1,\dots,X_k}(\mu)=(I(\pi_1)(\mu),\dots,I(\pi_k)(\mu)),\ \mu\in I(\prod X_i)$$
(here $\pi_j\colon \prod X_i\to X_j$ is the projection onto the
$j$th factor). It is proved in \cite{E1} that the corresponding
map is open for the case of the functor of probability measures.
The following simple example shows that this is no true for the
functor of idempotent probability measures.
\begin{exam} Let $X=\{x_1,x_2\}$, $Y=\{y_1,y_2\}$, and   $$\mu=0\odot\delta_{(x_1,y_1)}\oplus
0\odot\delta_{(x_2,y_2)}\in I(X\times Y).$$ Then
$$M_{X,Y}(\mu)=(\mu_1,\mu_2)=(0\odot\delta_{x_1}\oplus
0\odot\delta_{x_2},0\odot\delta_{y_1}\oplus 0\odot\delta_{y_2})\in
I(X)\times I(Y).$$ For every natural $l$, let
$$\mu_1^{(l)}=\left(-\frac1l\right)\odot\delta_{x_1}\oplus
0\odot\delta_{x_2},\ \mu_2^{(l)}=0\odot\delta_{y_1}\oplus
\left(-\frac1l\right)\odot\delta_{y_2}.$$ Then there is no
sequence $(\mu^{(l)})_{l=1}^\infty$ in $I(X\times Y)$ with
$M_{X,Y}(\mu^{(l)})=(\mu_1^{(l)},\mu_2^{(l)})$ and
$\lim_{l\to\infty}\mu^{(l)}=\mu$. This shows that the map
$M_{X,Y}$ is not open.

\end{exam}

A diagram (\ref{dia:1}) in the category $\comp$  is called {\em
open-bicommutative} if its characteristic map $\chi$ is an open
onto map. A functor $F$ acting in the category $\comp$ is called
{\em open-bicommutative} if it preserves the class of
open-bicommutative diagrams. See \cite{K} for the proof of
open-bicommutativity of some functors related to the
probability-measure functor.

Clearly, the product diagram
\begin{equation}\label{dia:2} \xymatrix{X\times Y\ar[r]\ar[d]&Y\ar[d]\\
 X\ar[r]&\{*\}}
 \end{equation}
 is open-bicommutative. It is clear from the above example that
 the diagram obtained by application the
 functor $I$ to diagram (\ref{dia:2}) is not open-bicommutative. This demonstrates that the
 functor $I$ is not open-bicommutative.
\section{Metrization}\label{s:metr}

Let $(X,d)$ be a compact metric space.

By $\mathrm{n-LIP}=\mathrm{n-LIP}(X,d)$ we denote the set of
Lipschitz functions with the Lipschitz constant $\le n$ from
$C(X)$.

Fix $n\in\mathbb N$. For every $\mu,\nu$, let $$\hat
d_n(\mu,\nu)=\sup\{|\mu(\varphi)-\nu(\varphi)|\mid \varphi\in
\mathrm{n-LIP}\}.$$

\begin{thm} The function $\hat d_n$ is a continuous pseudometric on $I(X)$.
\end{thm}
\begin{proof} We first remark that $\hat d_n$ is well-defined.
Indeed, $\sup\varphi-\inf\varphi\le n\diam X$, for every
$\varphi\in \mathrm{n-LIP}$, whence
$|\mu(\varphi)-\nu(\varphi)|\le2n\diam X$.

Obviously, $\hat d_n(\mu,\mu)=0$ and $\hat d_n(\mu,\nu)=\hat
d_n(\nu,\mu)$, for every $\mu,\nu\in I(X)$.

We are going to prove that $\hat d$ satisfies the triangle
inequality. Since, for every $\varphi\in \mathrm{n-LIP}\}$ and
$\mu,\nu,\tau\in I(X)$, $$\hat d(\mu,\nu)\ge
|\mu(\varphi)-\nu(\varphi)|,\ \hat d(\nu,\tau)\ge
|\nu(\varphi)-\tau(\varphi)|,$$ we have $$\hat d_n(\mu,\nu)+\hat
d(\nu,\tau)\ge|\mu(\varphi)-\nu(\varphi)|+
|\nu(\varphi)-\tau(\varphi)|\ge |\mu(\varphi)-\tau(\varphi)|,$$
whence, passing to $\sup$ in the right-hand side, we obtain $\hat
d_n(\mu,\nu)+\hat d(\nu,\tau)\ge\hat d(\mu,\tau)$.

Now, we prove that $\hat d$ is continuous. Suppose the contrary.
Then one can find a sequence $(\mu_i)_{i=1}^\infty$  in $I(X)$
such that $\lim_{i\to\infty}\mu_i=\mu\in I(X)$ and $\hat
d(\mu_i,\mu)\ge c'$, for some $c'>0$. Then there exist
$\varphi_i\in\mathrm{n-LIP}$, $i\in \mathbb N$,  such that
$|\mu_i(\varphi_i)-\mu(\varphi_i)|\ge c>$, for some $c>0$. Since
the functionals in $I(X)$ are weakly additive, without loss of
generality, one may assume that $\varphi_i(x_0)=0$, for some base
point $x_0\in X$, $i\in \mathbb N$. By the Arzela-Ascoli theorem,
there exists a limit point $\varphi\in\mathrm{n-LIP} $ of the
sequence  $(\varphi_i)_{i=1}^\infty$. We have
$|\mu_i(\varphi)-\mu(\varphi)|\ge c$, which contradicts to the
fact that $(\mu_i)_{i=1}^\infty$ converges to $\mu$.
\end{proof}
\begin{rem} Simple examples demonstrate that $\hat d$ cannot be a
metric whenever $X$  consists of more than one point.
\end{rem}

\begin{prop} The family of pseudometrics $\hat d_n$, $n\in\mathbb
N$, separates the points in $I(X)$.
\end{prop}

\begin{proof} Let $\mu,\nu\in I(X)$, $\mu\neq\nu$. There exists
$\varphi\in C(X)$ such that $|\mu(\varphi)-\nu(\varphi)|>c$, for
some $c>0$. There exists $\psi\in \mathrm{n-LIP}$, for some
$n\in\mathbb N$, such that $\|\varphi-\psi\|\le(c/3)$. Then,
clearly, $|\mu(\psi)-\nu(\psi)|\ge(c/3)$ and therefore $\hat
d_n(\mu,\nu)\ge(c/3)$.
\end{proof}

We let $\tilde d_n=(1/n)\hat d_n$.

\begin{prop} The map $\delta=\delta_X$, $x\mapsto\delta_x\colon
(X,d)\to(I(X),\tilde d_n)$, is an isometric embedding for every
$n\in\mathbb N$.
\end{prop}

\begin{proof} Let $x,y\in X$ and $\varphi\in \mathrm{n-LIP}$. Then
$|\delta_x(\varphi)-\delta_y(\varphi)|\le n d(x,y)$, therefore
$\hat d_n(\delta_x,\delta_y)\le n d(x,y)$. Thus $\tilde
d_n(\delta_x,\delta_y)\le  d(x,y)$.

On the other hand, define $\varphi_x\in\mathrm{n-LIP}$ by the
formula $\varphi_x(z)=nd(x,z)$, $z\in X$. Then
$|\delta_x(\varphi_x)-\delta_y(\varphi_x)|=nd(x,y)$ and we are
done.
\end{proof}

\begin{prop} Let $f\colon (X,d)\to (Y,\varrho)$ be a nonexpanding
map of compact metric spaces. Then the map $I(f)\colon (I(X),\hat
d_n)\to (I(Y),\hat\varrho_n)$ is also nonexpanding, for every
$n\in\mathbb N$.
\end{prop}
\begin{proof} Given $\varphi\in\mathrm{n-LIP}(Y)$, note that $\varphi
f\in\mathrm{n-LIP}(X)$ and, for any $\mu,\nu\in I(X)$, we have
$$|I(f)(\mu)(\varphi)-I(f)(\nu)(\varphi)|=|\mu(\varphi f)-\nu(\varphi f)|\le \hat d_n(\mu,\nu).$$
Passing to the limit in the left-hand side of the above formula,
we are done.
\end{proof}

Note that the above construction of $\hat d$ can be applied not
only to metrics but also to continuous pseudometrics.  Proceeding
in this way we obtain the iterations $(I(X),\tilde d_n)$,
$(I^2(X),\tilde{\tilde{d}}_{nm}=(\tilde{d}_n)\tilde{}_m)$,\dots

\begin{prop} For a metric space $(X,d)$, the map $\zeta_X\colon (I^2(X),\tilde{\tilde{d}}_{nn})
\to (I(X),\tilde d_n)$ is
nonexpanding.
\end{prop}
\begin{proof} We first prove that, for any $\varphi\in\mathrm{n-LIP}(X,d)$, we have
$\bar\varphi\in \mathrm{n-LIP}(I(X),\hat d)$. Indeed, given
$\mu,\nu\in I(X)$, we see that $$n\tilde d(\mu,\nu)= \hat
d(\mu,\nu)\ge
|\mu(\varphi)-\nu(\varphi)|=|\bar\varphi(\mu)-|\bar\varphi(\nu)|$$
and we are done.

Suppose now that $M,N\in I^2(X)$, $\mu=\zeta_X(M)$,
$\nu=\zeta_X(N)$. Given $\varphi\in\mathrm{n-LIP}(X,d)$, we obtain
$$|\mu(\varphi)-\nu(\varphi)|=|M(\bar\varphi)-N(\bar\varphi)|\le \tilde{\tilde{d}}_{nn}(M,N).$$
Passing to the limit in the left-hand side, we are done.

\end{proof}

\begin{rem} Using the results on existence of the pseudometrics
$\tilde d_n$, one can define the spaces of idempotent probability
measures with compact support for metric and, more generally,
uniform spaces. Indeed, let $(X,d)$ be a metric space. We define
the set $I(X)$ to be the direct limit of the direct system
$\{I(A),I(\iota_{AB});\exp X\}$ (here, for $A,B\in \exp X$ with
$A\subset B$, we denote by $\iota_{AB}\colon A\to B$ the inclusion
map). For every $A\in\exp X$, we identify $I(A)$ with the
corresponding subset of $I(X)$ along the map $I(\iota_A)$, where
$\iota_A\colon A\to X$ is the limit inclusion map. For any $\mu\in
I(X)$, there exists a unique minimal $A\in\exp X$ such that
$\mu\in I(A)$. Then we say that $A$ is the {\em support} of $\mu$
and write $\mathrm{supp}(\mu)=A$.

Now, define a family of pseudometrics $\hat d_n$, $n\in\mathbb N$,
on $I(X)$ as follows. Given $\mu,\nu\in I(X)$, we let $$\hat
d_n(\mu,\nu)= \hat
d_n|((\mathrm{supp}(\mu)\cup\mathrm{supp}(\nu))\times(\mathrm{supp}(\mu)\cup\mathrm{supp}(\nu)))(\mu,\nu).$$

One can prove that, for any uniform space $(X,\mathcal U)$, if the
uniformity $\mathcal U$ is generated by a family $\{d^\alpha\mid
\alpha\in A\}$ of pseudometrics, then the family $\{\tilde
d^\alpha_n\mid \alpha\in A,\ n\in\mathbb N\}$ of pseudometrics on
$I(X)$ generates a uniformity on $I(X)$.
\end{rem}

\section{Remarks and open problems}

L. Shapiro \cite{Sh} remarked that $P$ is the minimal normal
functor that admits a factorization through the category of
compact convex sets (in locally convex spaces) and affine
continuous maps.

\begin{que} Is $I$ the minimal normal functor that
 admits a factorization through the category of
compact max-plus convex sets?
\end{que}

\subsection{Idempotent barycentrically open max-plus convex sets}

V. Fedorchuk \cite{F2} characterized the barycentrically open
compact convex sets, i.e. compact convex sets $X$ for which the
barycenter map $P(X)\to X$ is open. Note that some
characterization results in this direction are also obtained in
\cite{E}, \cite{Pa}, \cite{OB}. In particular, it is proved in
\cite{OB} that a compact convex set $K$ in a locally convex space
is barycentrically open if the map $(x,y)\mapsto \frac12(x+y)$ is
open.

\begin{que} Characterize the class of max-plus convex compact
spaces for which the idempotent barycenter map is open. In
particular, is the latter property equivalent to the openness of
the map $(x,y)\mapsto x\oplus y$?
\end{que}

It is proved in \cite{F1} that the product of barycentrically open
compact convex sets is again barycentrically open.
\begin{que} Is an analogous fact true for idempotent barycentrically open max-plus convex
sets?
\end{que}
\subsection{Idempotent probability measure monad}
V. Fedorchuk \cite{F3} proved that there exists a unique monad in
$\comp$ with the probability measure functor as its functorial
part. It follows from the general properties of normal functors in
$\comp$ that there exists a unique natural transformation
$\id_{\comp}\to I$. This leads to the following question.
\begin{que} Is $\zeta\colon I^2\to I$ the unique natural
transformation that determines a monad structure for the functor
$I$?
\end{que}

T. \'Swirszcz \cite{Sw} proved that the category of compact convex
sets and affine continuous maps is monadic over the category
$\comp$. This leads to the following question.

\begin{que} Is the category of (suitably defined) compact max-plus convex
sets in locally convex lattices and affine continuous maps monadic
over the category $\comp$?
\end{que}

\begin{que} Characterize the category of $\mathbb I$-algebras.
\end{que}

\subsection{Milyutin maps}
We borrowed the idea of the proof of Theorem \ref{t:mil} from
\cite{AT}. Similarly as in \cite{AT}, one can prove that one can
choose a map $s\colon Y\to P(X)$  so that
$\mathrm{supp}(s(y))=f^{-1}(y)$, for any $y\in Y$, and, moreover,
every idempotent probability measure $s(y)$ is atomless in some
appropriate sense.

As we already remarked, it was first proved in \cite{DE} that the
probability measure functor $P$ preserves the class  of open maps.
The openness of the functor $O$ is proved in \cite{R}. The method
applied in \cite{DE} and \cite{R} does not work in our case.

The proof of Theorem \ref{t:open} is based on the properties of
Milyutin maps and can be also applied to the proof of openness of
the functor $P$ (see \cite{Z}) as well as of another related
functors.

Similarly like in \cite{SS}, one can use Milyutin maps in order to
prove a counterpart of the Michael selection theorem for the
max-plus-convex valued maps. That such a theorem can be proved by
methods based on general convexity structures is indicated in
\cite{BC} (see $\mathbb B$-spaces Metatheorem 5.0.19 therein).  We
return to this topic in another publication.

\subsection{Metrization}
If $(X,d)$ is a compact metric space, then the space $P(X)$ can be
endowed with the Kantorovich metric. It is an open problem whether
there exists a natural metrization of the space $I(X)$.

\begin{que} Is there a metrization of the space $I(X)$, for all
compact metric spaces $X$, which makes the monad $\mathbb I$
perfectly metrizable  (see \cite{F4} for the definition)?
\end{que}

Note that, for some natural reasons, to the notion of metric in
classical analysis and topology there correspond that of
ultrametric  in the idempotent case (see \cite{HZ}). (Recall that
a metric $d$ on a set $X$ is called an {\em ultrametric} if the
following strict triangle inequality holds:
$d(x,y)\le\max\{d(x,z),d(z,y)\}$). A counterpart of the
Kantorovich metric on the space $I(X)$ can be defined in the case
of an ultrametric space $X$. We leave it as an open problem to
define a Kantorovich-type metric on the spaces of idempotent
probability measures. V. Fedorchuk \cite{F4} introduced the notion
of perfectly metrizable monad. Roughly speaking, this is a monad
$(F,\eta,\mu)$ on $\comp$ whose functorial part $F$ is metrizable
(i.e., it can be lifted to the category of compact metric  spaces
and nonexpanding maps) and the maps $\eta_X$ and $\mu_X$ are
nonexpanding. The results of Section \ref{s:metr} suggest that one
can introduce an analogous structure of monad metrizable by a
countable family of pseudometrics.

\begin{que} Is there a counterpart of the Prokhorov metric for the
functor of idempotent probability measures?
\end{que}

\subsection{Idempotent probability measures of noncompact spaces}
Let $X$ be a Tychonov space and $\beta X$ be its Stone-\v Cech
compactification. We define $I(X)=\{\mu\in I(\beta X)\mid
\mathrm{supp}(\mu)\subset X\}$. For any maps $f\colon X\to Y$ of
Tychonov spaces, we have $I(\beta f)(I(X))\subset I(Y)$ and,
therefore, we can define the map $I(f)= I(\beta f)|I(X)\colon
I(X)\to I(Y)$. We thus obtain an extension of $I$ onto the
category of Tychonov spaces. Note that a base of topology on
$I(X)$ can be formed by the sets of the form $$\langle\mu;
\varphi_1,\dots,\varphi_n; \varepsilon\rangle=\{\nu\in I(X)\mid
|\mu(\varphi_i)-\nu(\varphi_i)|<\varepsilon,\ i=1,\dots,n\},$$
where  $\varphi_i\in C(X)$, $i=1,\dots,n$, are bounded and
$\varepsilon>0$.

\subsection{Generalizations} One can obtain counterparts of the
above results for another spaces of pseudo-additive measures. Here
we mention only one example. Let $\cdot$ denote the $\min$
operation on the set $[-\infty,\infty]$. The set of functionals
$\mu\colon C(X)\to\mathbb R$ satisfying
$\mu(\varphi)\oplus\mu(\psi)$, $\mu(c_X)=c$, and
$\mu(\lambda\cdot\varphi)=\lambda\cdot\mu(\varphi)$ can be
topologized by the weak* topology and some of the results of this
paper have their counterparts for such functionals.

\subsection{Geometric properties}   Many publications are devoted
to geometric properties of the probability measure functor. Some
of them have their counterparts for the idempotent probability
measures. We will consider some these properties in subsequent
publications.

\end{document}